\documentclass[a4paper,reqno,10pt,oneside]{amsart}
\usepackage{amsmath,geometry}
\usepackage{graphicx}
\usepackage{mathrsfs}
\usepackage{amssymb,amsmath, amsfonts, amsthm}
\usepackage{latexsym}
\usepackage{color} 
\usepackage[dvips]{epsfig}
\usepackage{graphicx}

\allowdisplaybreaks[1]

\numberwithin{equation}{section}

\renewcommand{\(}{\left(}
\renewcommand{\)}{\right)}
\renewcommand{\[}{\left[}
\renewcommand{\]}{\right]}

\newtheorem{theorem}{Theorem}[section]
\newtheorem{proposition}[theorem]{Proposition}
\newtheorem{corollary}[theorem]{Corollary}

\newtheorem{remark}[theorem]{Remark}
\newtheorem{definition}[theorem]{Definition}

\renewcommand{\le}{\leqslant}
\renewcommand{\ge}{\geqslant}

\newcommand{\beq}{\begin{equation}}
\newcommand{\eeq}{\end{equation}}
\newcommand{\beqs}{\begin{equation*}}
\newcommand{\eeqs}{\end{equation*}}
\newcommand{\beqn}{\begin{eqnarray}}
\newcommand{\eeqn}{\end{eqnarray}}
\newcommand{\beqns}{\begin{eqnarray*}}
\newcommand{\eeqns}{\end{eqnarray*}}
\newcommand{\bdoc}{\begin{document}}
\newcommand{\edoc}{\end{document}}
\newcommand{\be}{\begin{enumerate}}
\newcommand{\ee}{\end{enumerate}}
\newcommand{\bdescr}{\begin{description}}
\newcommand{\edescr}{\end{description}}
\newcommand{\ba}{\begin{array}}
\newcommand{\ea}{\end{array}}
\newcommand{\intR}{\int_{\mathbb R^N}}
\newcommand{\R}{\mathbb R}
\newcommand{\C}{\mathbb C}
\newcommand{\parallelsum}{\mathbin{\!/\mkern-5mu/\!}}
\newcommand{\e}{\epsilon}
 \renewcommand{\(}{\left(}
\renewcommand{\)}{\right)}
\renewcommand{\[}{\left[}
\renewcommand{\]}{\right]}
\newenvironment{Proof}{\noindent{\bf Proof}}{\hfill$\Box$\\[2mm]}

\begin{document}
\title[H-surfaces in cones]{Existence of stable H-surfaces in cones and\\ their representation as radial graphs}

\author{P\MakeLowercase{aolo} Caldiroli, A\MakeLowercase{lessandro} Iacopetti}

\subjclass[2010]{53A10, 35J66, 57R40}
\keywords{Prescribed mean curvature, Plateau's problem, H-surfaces}
\thanks{\emph{Acknowledgements.} Research partially supported by the project ERC Advanced Grant 2013 n.~339958 Complex Patterns for Strongly Interacting Dynamical Systems COMPAT, and by Gruppo Nazionale per l'Analisi Matematica, la Pro\-ba\-bi\-li\-t\`a e le loro Applicazioni (GNAMPA) of the Istituto Nazionale di Alta Matematica (INdAM). The first author is also supported by the PRIN-2012-74FYK7 Grant ``Variational and perturbative aspects of nonlinear differential problems''}
\address[Paolo Caldiroli]{Dipartimento di Matematica, Universit\`a di Torino, via Carlo Alberto, 10 -- 10123 Torino, Italy}
\email{paolo.caldiroli@unito.it}
\address[Alessandro Iacopetti]{Dipartimento di Matematica, Universit\`a di Torino, via Carlo Alberto, 10 -- 10123 Torino, Italy}
\email{alessandro.iacopetti@unito.it}
\begin{abstract}
In this paper we study the Plateau problem for disk-type surfaces contained in conic regions of $\R^{3}$ and with prescribed mean curvature $H$. Assuming a suitable growth condition on $H$, we prove existence of a least energy $H$-surface $X$ spanning an arbitrary Jordan curve $\Gamma$ taken in the cone. Then we address the problem of describing such surface $X$ as radial graph when the Jordan curve $\Gamma$ admits a radial representation. Assuming a suitable monotonicity condition on the mapping $\lambda\mapsto \lambda H(\lambda p)$ and some strong convexity-type condition on the radial projection of the Jordan curve $\Gamma$, we show that the $H$-surface $X$ can be represented as a radial graph.  
\end{abstract}

\maketitle

\section{Introduction}

In the present paper we aim to investigate some aspects on the Plateau problem for disk-type surfaces with prescribed mean curvature in the directions described as follows. 
%
%
%
Fixing a cone of angular radius $\beta$
$$
\mathfrak{C}_{\beta}:=\{p=(x,y,z)\in\R^3~|~ z>|p| \cos \beta\},
$$ 
a Jordan curve $\Gamma\subset\overline{\mathfrak{C}_{\beta}}\setminus\{0\}$, and a mapping $H\colon\overline{{\mathfrak{C}_{\beta}}}\to\R$, we are interested in finding conditions on $H$, possibly related to $\beta$, ensuring that stable surfaces in $\overline{\mathfrak{C}_{\beta}}\setminus\{0\}$ with mean curvature $H$, spanning $\lambda\Gamma$ do exist for every $\lambda>0$. Moreover we address the problem of describing such surfaces as radial graphs when their boundaries admit a radial representation. 

In order to state our main results, let us state the analytical formulation of the problem. Let $B=\{(u,v)\in\R^{2}~|~u^{2}+v^{2}<1\}$ be the unit open disk. In general, the Plateau problem for a given Jordan curve $\Gamma$ and a prescribed mean curvature function $H$ consists in looking 
for maps $X\colon\overline{B}\to\R^{3}$ solving
\begin{align}
\label{eqprob}
&\Delta X = 2 H(X) X_u\wedge X_v~~\text{ in $B$}\\
\label{conf}
&|X_u|^2-|X_v|^2=0=X_u \cdot X_v~~\text{ in $B$}\\
\label{eq:boundary}
&X|_{\partial B}\colon\partial B\to\Gamma\text{ is an (oriented) parametrization of }\Gamma.
\end{align}
A map $X \in C^0(\overline{B},\R^3)\cap C^2(B,\R^3)$ satisfying (\ref{eqprob})--(\ref{eq:boundary}) will be called \emph{$H$-surface spanning $\Gamma$} (see \cite{Struwe}). It is known that that if $X$ is an $H$-surface, then $X$ has mean curvature $H(X)$ apart from branch points, i.e., points $(u,v)\in B$ where $\nabla X(u,v)=0$. 

Our first result can be stated as follows.

\begin{theorem}
\label{thm:existence}
Let $\beta\in(0,\frac{\pi}{2})$ and let $H\colon\overline{\mathfrak{C}_{\beta}}\to\R$ be a mapping of class $C^{1}$, satisfying
\begin{equation}
\label{growthcondition}
|H(p)||p|\le\frac{\cos\beta}{2(1+\cos\beta)}\quad\forall p\in\overline{\mathfrak{C}_{\beta}}~\!.
\end{equation}
Then for every rectifiable Jordan curve $\Gamma\subset\overline{\mathfrak{C}_{\beta}}\setminus\{0\}$ there exists an $H$-surface $X \in C^0(\overline{B},\R^3) \cap C^2(B,\R^3)$ spanning $\Gamma$ and contained in $\overline{\mathfrak{C}_{\beta}}\setminus\{0\}$. Moreover we have that $X(B) \subset \mathfrak{C}_\beta$.  
\end{theorem}
We point out that the assumption (\ref{growthcondition}) fixes a bound on the radial behaviour of $H$ with respect to the angular diameter of the given Jordan curve $\Gamma$. Moreover, since (\ref{growthcondition}) is asked to hold on a dilation-invariant domain and is independent of the curve $\Gamma$, the existence result stated by Theorem \ref{thm:existence} remains true also taking $\lambda\Gamma$ instead of $\Gamma$, for every $\lambda>0$. Note that the case of nonzero constant mean curvature is ruled-out. 

In fact we can provide more information on the $H$-surface given by Theorem \ref{thm:existence}. More precisely, taking the variational character of the Plateau problem into account, such $H$-surface is  characterized as a least energy surface, namely is a minimum point of the energy functional 
associated to system (\ref{eqprob}), in the class of admissible mappings 
satisfying (\ref{eq:boundary}). We refer to Sections 2 and 3 for more details about this aspect.

Our second result provides an answer to the issue of representing an $H$-surface as a radial graph, when its contour is a radial graph. To this purpose, we need a monotonicity condition on the mapping $\lambda\mapsto \lambda H(\lambda p)$ and some strong convexity-type condition on the radial projection of the Jordan curve $\Gamma$. In particular, we can show:

\begin{theorem}
\label{thm:radial-graph}
Let $\beta\in(0,\frac{\pi}{2})$ and let $H\colon\overline{\mathfrak{C}_{\beta}}\to\R$ be a mapping of class $C^{1,\alpha}$, satisfying (\ref{growthcondition}) and 
\begin{equation}
\label{monassumption}
H(p)+\nabla H(p)\cdot p\ge 0\quad\forall p\in\overline{{\mathfrak{C}}_{\beta}}~\!.
\end{equation}
Let $\Gamma$ be a regular Jordan curve of class $C^{3,\alpha}$ contained in $\overline{\mathfrak{C}_{\beta}}\setminus\{0\}$ and let $X$ be the least energy $H$-surface spanning $\Gamma$, given by Theorem \ref{thm:existence}. Assume that:
\begin{itemize}
\item[(i)] $\Gamma$ is a radial graph, i.e. there exists a domain $\Omega\subset\mathbb{S}^{2}$ and a map $g\colon\partial\Omega\to\R^{+}$ (with the same regularity of $\Gamma$) such that $\Gamma=\{g(p)p~|~p\in\partial\Omega\}$;
\item[(ii)] the domain $\Omega$ is $\beta$-convex (see Definition \ref{defbetaconvx});
\item[(iii)] the radial projection of $X|_{\partial B}$ induces a positive orientation on $\partial\Omega$ (see Definition \ref{definitposorjc} and  Remark \ref{remarkposorientms}).
\end{itemize} 
Then the radial projection of $X$ is a diffeomorphism between $\overline{B}$ and $\overline\Omega$ and $X(\overline{B})$ can be represented as a radial graph. In particular $X$ has no branch point.  
\end{theorem}
\noindent
We notice that Theorem \ref{thm:radial-graph} is a corollary of a more general result (Theorem \ref{mainteoglobinv}) about the representation of stable $H$-surfaces as radial graphs. The meaning of \emph{stable} $H$-surface is explained in Definition~\ref{defstablehsurf}.

The study developed in the present paper is the natural counterpart of analogous issues on the Plateau problem for disk-type $H$-surfaces in a cylinder and their representation as cartesian graphs with respect to the direction of the axis of the cylinder. On this side some results are already known in in the literature: Rad\'o proved in \cite{Rado} that minimal surfaces, spanning a Jordan curve with one-one projection onto the boundary of a planar convex domain $D \subset \R^2$, can be represented as cartesian graphs of a function over $D$. Serrin in \cite{Serrin1}, Gulliver and Spruck in \cite{GulliverSpruck1} proved the same result in the case of surfaces of constant mean curvature, but with different assumptions. Sauvigny in \cite{Sauvigny82} studied the case of stable $H$-surfaces with $H$ not necessarily constant. In particular, he proved that, under a suitably strong convexity condition on the planar domain $D$ (which is the planar version of our $\beta$-convexity property), if $H$ is monotone along the axial direction of the cylinder, then stable $H$-surfaces can be represented as cartesian graphs of a real valued function over $D$.

Clearly our results cannot be recovered by those obtained for the cylinder. Indeed we deal with a problem whose geometry exhibits a dilation invariance, in the sense that conditions (\ref{growthcondition}) and (\ref{monassumption}) regard just the radial behaviour of the prescribed mean curvature function $H$. 

Furthermore, considering cones and radial projections rather that cylinders and cartesian projections lead some non trivial, extra difficulties. The reason is that conical surfaces exhibit a singular point at their vertex and the radial projection is a nonlinear mapping. Let us display the main difficulties by highlighting the more delicate steps in our arguments. 

Concerning the existence result stated by Theorem \ref{thm:existence}, we follow the standard procedure of minimizing the energy functional associated to the $H$-system (\ref{eqprob}) in the class of admissible functions. Assumption (\ref{growthcondition}) guarantees that the energy functional is bounded from below and, by known results, one gets existence of a minimizer $X$. Actually, in principle, the minimizer could touch the obstacle, in particular the vertex of the cone $\mathfrak{C}_{\beta}$. To overcome this difficulty we smooth the cone at the origin in a suitable way and we use a deep result by Gulliver and Spruck (see \cite{GulliverSpruck2}), together with the growth condition \eqref{growthcondition}, in order to obtain that the minimizer does not touch the boundary of the smoothed cone and then stays far from the vertex of $\mathfrak{C}_{\beta}$. Thus, by well known regularity results (see for instance \cite{RegMinSurf}), $X$ is a classical solution of \eqref{eqprob}--\eqref{eq:boundary}. Notice that our procedure needs more care than in the case of the analogous obstacle problem in a ball or in a cylinder (see Theorems 8 and 9 in Section 4.7 of \cite{RegMinSurf}). We also observe that the minimizer $X$ turns out to be stable in the sense of Definition \ref{defstablehsurf}, provided that $\Gamma$ and $H$ are regular enough (see also Proposition \ref{regpropbhdata}).
 
Now let us spend a few words about our second result, concerning the characterization of a stable $H$-surface $X$ as radial graph, and let us shortly illustrate the strategy followed to show that the radial projection $PX={X}/{|X|}$ is a homeomorphism between $\overline B$ and $\overline \Omega=PX(\overline{B})$. Under the assumptions on $H$ and $\Gamma$ as in the statement of Theorem \ref{thm:radial-graph}, we prove that the radial component of the Gauss Map $N$ is always positive in $\overline B$, namely
\beq\label{nscalx}
N \cdot X>0 \ \ \ \hbox{in} \ \overline B. 
\eeq
The maximum principle is the key tool to this aim. In fact, property \eqref{nscalx} implies local invertibility of $PX$ far from branch points. The issue of global invertibility is not tackled with the same strategy followed for the analogous problem of the projection along a fixed direction as in the papers \cite{GulliverSpruck1} and \cite{Sauvigny82}, because the expansion about branch points, based on the Hartman-Wintener technique, does not fit well with the radial projection. Instead, we follow an argument which is mainly based on the degree theory, combined with a classical result about global invertibility (see \cite{AmbrosettiProdi}) and Jordan-Sch\"onflies's Theorem (see \cite{Thomassen}).

Finally we notice that the (non-parametric) Plateau problem for $H$-surfaces characterized as radial graphs was already discussed by Serrin in in \cite{Serrin2}. Actually in that work a class of positively homogeneous prescribed mean curvature functions is considered and the existence of $(n-1)$-dimensional $H$-surfaces in $\R^n$ spanning a datum $\Gamma$ is proved under the following assumptions: $\Gamma$ is the radial graph of a positive mapping $f$ defined on the boundary of a given smooth domain $\Omega$ contained in a hemisphere of $\mathbb{S}^n$, and  
$$
H_g(y) \geq \frac{n}{n-1} H(y) f(y)\quad\forall y \in \partial \Omega~\!,
$$
where $H_g$ denotes the geodesic mean curvature of $\partial \Omega$.
For spherical caps, this condition turns out to be less restrictive than \eqref{growthcondition}. On the other hand, our results allow sign-changing and non-homogeneous mean curvature functions, which cannot be considered in \cite{Serrin2}.

Lastly, let us sketch an outline of the present paper: Sect.~2 contains a collection of known facts and technical results which will be used in the sequel. In Sect.~3 we prove the existence result stated by Theorem \ref{thm:existence}. In Sect.~4 we discuss the notion of $\beta$-convex domain in $\mathbb{S}^{2}$. Finally Sections 5 and 6 contain the proof of \eqref{nscalx} and of Theorem \ref{thm:radial-graph} (actually, a more general version), respectively.

\section{Notation and preliminary results}
In this section we fix some notation and we collect some known facts which will be useful in the rest of the paper.

We denote by $B$ the unit open disk of $\R^2$ and by $\overline{B}$ its closure. 
 We will use indistinctly both the real notation $(u,v)$ or the complex notation $z$, $w$ to denote a generic point of $B$ or $\overline{B}$. In particular, it will be always understood that $z=e^{i \theta} \in \partial B$ stands for $(\cos \theta, \sin \theta)$. We denote by $\mathbb{S}^2$ the unit sphere of $\R^3$ and by $P\colon\R^3\setminus\{0\} \to \mathbb{S}^2$ the radial projection map, defined by $P(x):=\frac{x}{|x|}$. We will use both the notation $P(Y)$ or $PY$ to denote the composition $P\circ Y$, whenever $Y$ is map with values in $\R^3\setminus\{0\}$.\\

We begin with recalling some important facts about branch points and the normal $N$ to an H-surface.
\begin{theorem}[see Theorem 1, Sect. 2.10, \cite{RegMinSurf} and also Remark 3, Sect. 5.1, \cite{MinSurf}] \label{HWexpansion}
 Let $X$ be an H-surface of class $C^{2,\alpha}(B, \R^3)$ or $C^{2,\alpha}(\overline{B}, \R^3)$ respectively. Then, for each point $w_0 \in B$ or $\overline{B}$, there is a vector $A=(A_1,A_2,A_3) \in \mathbb{C}^3$ with $A_1^2+A_2^2+A_3^2=0$, and a nonnegative integer $n=n(w_0)$ such that

\begin{equation}\label{esphw}
X_w(w)=A(w-w_0)^n + o(|w-w_0|^n), \ \ \hbox{as} \ w \to w_0,
\eeq
where $X_w:=\frac{1}{2}(X_u-iX_v)$.
\end{theorem}

\begin{remark}
The point $w_0$ in the above statement is a branch point of $X$ if and only if $n(w_0)\geq 1$, and in this case $n(w_0)$ is called the order of the branch point $w_0 \in B$ (or $\overline{B}$ respectively). Obviously $w_0$ is regular point of $X$ if and only if $n(w_0)=0$. Thanks to \eqref{esphw} we deduce that branch points of an H-surface are isolated. In particular, if $X \in C^{2,\alpha}(\overline{B}, \R^3)$ then the set of branch points is finite. 
\end{remark}

In order to get more analytic regularity on the solution $X$ we have to ask more regularity on the function $H$ and on the Jordan curve $\Gamma$.  More precisely, we recall that:
\begin{proposition}[see Chap. IX, Sect. 4, \cite{Sauvigny2006} and Sect. 2.3, \cite{RegMinSurf}] \label{regpropbhdata}
\ \ \ \ \ \ \ 
\begin{itemize}
\item[(i)]If $H \in C^{r,\alpha}(\R^3)$, for $r \in \mathbb{N}$, $\alpha \in (0,1)$, then any solution $X \in C^2(B,\R^3)$ of \eqref{eqprob} is of class $C^{r+2,\alpha}(B,\R^3)$.
\item[(ii)] If $H \in C^{0,\alpha}(\R^3)$, for some $\alpha \in (0,1)$, and $X$ is an H-surface such that $X(\partial B)$ lies on a regular Jordan curve of class $C^{2,\alpha}$ then $X \in C^{2,\alpha}(\overline{B}, \R^3)$. More in general, if $H \in C^{r-2,\alpha}(\R^3)$ and $\Gamma \in C^{r,\alpha}$, for some $r\geq 2$, then $X \in C^{r,\alpha}(\overline{B}, \R^3)$. 
\end{itemize}
 \end{proposition}

For $X \in C^{2,\alpha}(B, \R^3)$ we denote by ${B}^\prime$ the set of regular points. We recall that for an H-surface $X \in C^{2,\alpha}(B, \R^3)$ and $w \in {B}^\prime$ the normal $N$ at $w$ is given by
\beq\label{normal}
N(w)= \frac{X_u(w)\wedge X_v(w)}{|X_u(w)\wedge X_v(w)|}=\frac{X_u(w)\wedge X_v(w)}{|X_u(w)|^2}.
\eeq

Thanks to the expansion \eqref{esphw}, writing $A=a - ib$, with $a, b \in \R^3$, from $A\neq 0$ and being $A_1^2+A_2^2+A_3^2=0$ it follows that $|a|=|b|\neq 0$, $a \cdot b=0$. Hence, if $w_0$ is a branch point, then

$$ N(w) \to \frac{a\wedge b}{|a|^2} \in \mathbb{S}^2, \ \ \hbox{as} \ w \to w_0, \ w \in {B}^\prime.$$

Therefore we deduce that the normal $N$ can be extended to a continuous function $N \in C^0(\overline{B}, \R^3)$ with $N(\overline{B}) \subset \mathbb{S}^2$. Furthermore we have:

\begin{theorem}[see Theorem 1, Sect. 5.1, \cite{MinSurf}]\label{teoeqnorm}
Assume that $H \in C^{1,\alpha}(\R^3)$ and that $X$ is an H-surface of class $C^{3,\alpha}(B,\R^3)$. Then the normal $N$  is of class $C^{2,\alpha}(B,\R^3)$ and satisfies the differential equation
\beq\label{eqsoddnorm}
\Delta N + 2 pN = - 2E \nabla H (X), 
\eeq
where $E:=|X_u|^2$,
\beq\label{densityfunct}
 p:= E \ [2H^2(X) - K -  (\nabla H(X)\cdot N)],
\eeq
is the so-called ``density function'' associated to $X$ and $K$ is the Gaussian curvature of $X$. Moreover $p \in C^{0,\alpha}(B)$.
\end{theorem}

As remarked in the introduction, it is well known that H-surfaces are obtained as stationary points of the energy functional 
\beq\label{enfunct}
\mathcal{F}(X)= \frac{1}{2}\int_B |\nabla X|^2 \ du \ dv + 2\int_B Q(X) \cdot X_u \wedge X_v \ du \ dv,
\eeq
where $Q\colon\R^3 \to \R^3$ is a vector field such that $div\ Q = H$. Let us also introduce the functional 
$$\mathcal{G}(X):=\int_B |X_u \wedge X_v|^2 \ du \ dv + 2\int_B Q(X) \cdot X_u \wedge X_v \ du \ dv.$$ Obviously we have $\mathcal{F}(X) \leq \mathcal{G}(X)$ and the equality $\mathcal{F}(X)=\mathcal{G}(X)$ holds if and only $X$ satisfies the conformality relations \eqref{conf}.\\

\begin{definition}
Let $X$ an H-surface of class $C^{3,\alpha}(\overline{B},\R^3)$ and let $\varphi \in C_0^\infty(B)$ be a test function.
We define the normal variation as the function $Z: \overline{B} \times (-\e_0,\e_0) \to \R^3$, $\e_0>0$, defined by 
$$ Z(w,t):=X(w) + t\varphi(w) N(w).$$
We define the first variation and the second variation of $\mathcal{G}$, in the normal direction $Z$, respectively, as
$$  \delta \mathcal{G}(X,\varphi N):=\frac{d}{dt}\mathcal{G}(Z) \Big{|}_{t=0}\quad\text{and}\quad
\delta^2 \mathcal{G}(X,\varphi N):=\frac{d^2}{dt^2}\mathcal{G}(Z) \Big{|}_{t=0}~\!.
$$
\end{definition}
The following result holds:
\begin{theorem}[see Theorem 1, Sect. 5.3, \cite{MinSurf}] \label{teonormvar}
Let $X \in C^{3,\alpha}(\overline{B},\R^3)$ an H-surface and let $\varphi \in C_0^\infty(B)$ a test function. Then: 
\begin{itemize}
\item[(i)] $\delta \mathcal{G}(X,\varphi N)=0$,
\item[(ii)] $\displaystyle \delta^2 \mathcal{G}(X,\varphi N)= \int_B |\nabla \varphi|^2 - 2 p \varphi^2 \ du \ dv$, where $p:B \rightarrow \R$ is given by \eqref{densityfunct}.
\end{itemize}
\end{theorem}
%

We recall now the fundamental notion of stability for H-surfaces.
\begin{definition}\label{defstablehsurf}
We say that an H-surface $X \in C^{3,\alpha}(\overline{B},\R^3)$ is stable if it satisfies the following inequality
$$ \delta^2 \mathcal{G}(X,\varphi N) \geq 0, \ \ \hbox{for all} \ \varphi \in C_0^{\infty}(B),$$
which, in view of Theorem \ref{teonormvar}, can be rewritten as
$$\displaystyle \int_B |\nabla \varphi|^2 - 2 p \varphi^2 \ du \ dv \geq 0, \ \ \hbox{for all} \ \varphi \in C_0^{\infty}(B).$$
\end{definition}

\begin{remark}\label{stabilitymin}
We point out that global and local minimizers of $\mathcal{G}$ are stable. In particular if $X$ satisfies the conformality relations and it is a minimizer of $\mathcal{F}$, then it is stable.
\end{remark}

The following result is a well known version of the maximum principle.
\begin{proposition}[see Proposition 1, Sect. 5.3, \cite{MinSurf}] \label{propstabsegnpos} 
Assume that $q \in C^{0,\alpha}(B)$ satisfies the stability inequality
$$\displaystyle \int_B |\nabla \varphi|^2 - 2 q \varphi^2 \ du \ dv \geq 0, \ \ \hbox{for all} \ \varphi \in C_0^{\infty}(B),$$
and let $f \in C^0(\overline{B})\cap C^2(B)$ be a solution of the boundary value problem
\beq\label{probstab}
\begin{cases}
\Delta f + 2qf \leq 0 & \hbox{in} \ B\\
f(w)>0 &  \hbox{on} \ \partial B.
\end{cases}
\eeq
Then $f(w)>0$ for all $w \in \overline{B}$.
\end{proposition}


Finally we state a classical result of global invertibility. We recall that a map $F\colon\mathcal{X}\to \mathcal{Y}$ between two topological spaces $\mathcal{X}$, $\mathcal{Y}$ is said to be proper if $F^{-1}(\mathcal{K})$ is compact in $\mathcal{X}$, for any compact subset $\mathcal{K}$ in $\mathcal{Y}$.

\begin{theorem}[see Theorem 1.8, Sect. 3.1, \cite{AmbrosettiProdi}] \label{globalinvtheo}
Let $\mathcal{X}$, $\mathcal{Y}$ be two Banach spaces and let $F\colon\mathcal{X}\to \mathcal{Y}$ be a continuous surjective proper map. Suppose that $F$ is locally invertible, $\mathcal{X}$ arcwise connected and $\mathcal{Y}$ simply connected. Then $F$ is a homeomorphism.
\end{theorem}

\section{Existence of $H$-surfaces in cones}

In this section we prove Theorem \ref{thm:existence}.  
We divide the proof in several steps.\\

\noindent\textbf{Step 1}: Extension of $H$ to $\overline{\mathfrak{C}_{\beta+ \delta}}$, for $\delta>0$ sufficiently small.\\
Let $\beta \in (0,\pi/2)$ and set 
$$
c_\beta:=\frac{\cos \beta}{2(1+\cos \beta)}~\!.
$$
Let $\bar\delta>0$ sufficiently small so that $\beta \pm \bar\delta \in (0,\pi/2)$ and $ c_{\beta-\delta} < \frac{\mathrm{cotan}{(\beta+\delta})}{2}$ for all $0<\delta<\bar\delta$ (this choice will be useful in the sequel of the proof). We point out that there always exists a $\bar\delta=\bar\delta(\beta)>0$ satisfying the previous inequality: in fact observe that since $\beta \in (0,\pi/2)$ the inequality $\frac{\cos \beta}{2(1+\cos \beta)} < \frac{\mathrm{cotan}{(\beta})}{2}$ is equivalent to $\sin \beta  < 1 + \cos \beta$, which holds true. Thus, by continuity of the function $\delta \mapsto \frac{\cos (\beta-\delta)}{2(1+\cos (\beta-\delta))} - \frac{\mathrm{cotan}{(\beta + \delta)}}{2}$ at $0$ we get the desired assertion.\\

 Let $0<\delta<\bar\delta$ sufficiently small so that $H$ can be extended to a function $H \in C^{1}( \overline{{\mathfrak{C}_{\beta+ \delta}}})$ with $|H(p)||p|\leq c_{\beta - \delta}$ (we observe that being $\gamma \mapsto c_\gamma$ a strictly decreasing function it holds $c_{\beta-\delta}>c_\beta$). Clearly $\Gamma$ is strictly contained in ${{\mathfrak{C}_{\beta+ \delta}}}$.\\

\noindent\textbf{Step 2:} Construction of a suitable smooth surface of revolution which approximates $\partial\mathfrak{C}_{\beta+\delta}$.\\
 The cone  $\partial\mathfrak{C}_{\beta+\delta}$ is a non-smooth surface of revolution obtained by rotating the half-line $\sigma(t)=(\sin(\beta+\delta) \ t,0,\cos (\beta+\delta) \ t)$, $t \in \R^+ \cup \{0\}$, lying in the $xz$-plane, through the $z$-axis.
 We consider the following approximating surface of revolution: let $\epsilon >0$ be a small parameter to be chosen later, 
 and set $t_\epsilon:= \frac{8}{3\sqrt{3}}\frac{1}{\cos(\beta+\delta)} \epsilon$, let $S_{\beta + \delta,\epsilon}$ be the surface obtained by rotating the curve $\sigma_\e(t):=(\alpha_1(t),0, \alpha_2(t))$ through the $z$-axis, parametrized by $\phi(t,\theta)=(\alpha_1(t)\cos \theta, \alpha_1(t) \sin \theta, \alpha_2(t))$, where
\beq\label{alpha12}
\begin{array}{lll}
 \displaystyle \alpha_1(t)&:=& \displaystyle  \sin(\beta+\delta) t, \ \ \ t \geq 0,\\[6pt]
 \displaystyle \alpha_2(t)&:=&\begin{cases}
 \displaystyle a_\e t ^4 + b_\e t ^2 + c_\e, &  \hbox{if} \ t \in [0, t_\epsilon], \\
 \cos (\beta+\delta) \ t, & \hbox{if} \ t \in]t_\epsilon , +\infty[. 
 \end{cases} 
\end{array}
\eeq
with $a_\e$, $b_\e$, $c_\e$ chosen in a suitable way in order that:
\begin{itemize}
\item[(i)] $S_{\beta + \delta,\epsilon}$ is of class $C^2$,
\item[(ii)] $0 \not\in S_{\beta + \delta,\epsilon}$,
\item[(iii)] the component $\mathfrak{S}_{\beta+\delta, \epsilon}$ of $\R^3\setminus {S}_{\beta+\delta, \epsilon}$ which does not contain the origin is convex,
\item[(iv)] the mean curvature (with respect to the inward normal) $H_{S_{\beta + \delta,\epsilon}}$ of $S_{\beta + \delta,\epsilon}$ satisfies
\beq\label{eqstimaH}
|H(p)| < H_{S_{\beta + \delta,\epsilon}}(p), \ \ \ \hbox{for any}\ p \in S_{\beta + \delta,\epsilon}.
\eeq
\end{itemize}

A good choice of the coefficients $a_\e$, $b_\e$, $c_\e$ is
$$ a_\e:=- \sqrt{3} \left(\frac{3}{8}\right)^4 \frac{ \cos^4(\beta+\delta)}{\epsilon^3}, \ b_\e:=2\sqrt{3} \left(\frac{3}{8}\right)^2  \frac{ \cos^2(\beta+\delta)}{\epsilon}, \ c_\e:= \frac{\epsilon}{\sqrt{3}}.  $$
Inequality \eqref{eqstimaH} is checked at Step 8.\\

\noindent\textbf{Step 3:} Choice of a vector field $Q: \overline{{\mathfrak{C}_{\beta+\delta}}} \to \R^3$ such that $div\ Q=H$.\\
Let us set $$Q(p):= \left(\int_{0}^1 H(tp)t^2 \ dt\right) \ p,\ \ \ p \in \overline{{\mathfrak{C}_{\beta+\delta}}}.$$
It is clear that $Q \in C^{1}(\overline{{\mathfrak{C}_{\beta+\delta}}},\R^3)$ and by elementary computations we see that $ div \ Q =  H$ in $\overline{{\mathfrak{C}_{\beta+\delta}}}$. 
Moreover we observe that, since $|H(p)| |p| \leq c_{\beta-\delta}$ for all $p \in \overline{{\mathfrak{C}_{\beta+\delta}}}$, we have that
\begin{equation}\label{eq1normlinftyQ}
\|Q\|_{\infty,  \overline{{\mathfrak{C}_{\beta+\delta}}}} \leq  \frac{ c_{\beta-\delta}}{2} < \frac{1}{4}.\\ 
\end{equation}

\noindent\textbf{Step 4:} Construction of a weak solution of \eqref{eqprob} which satisfies \eqref{conf} a.e. in $B$.\\
Let $\e>0$ sufficiently small such that $\Gamma \subset {\mathfrak{S}_{\beta + \delta,\epsilon}}$ and \eqref{eqstimaH} holds. We consider the variational problem $\mathcal{P}(\Gamma,\overline{\mathfrak{S}_{\beta,\epsilon}})$ given by $$\min_{X \in \mathcal{C}(\Gamma,\overline{\mathfrak{S}_{\beta + \delta,\epsilon}})} \mathcal{F}(X),$$
where $ \mathcal{F}$ is the functional defined in \eqref{enfunct} and $\mathcal{C}(\Gamma,\overline{\mathfrak{S}_{\beta+\delta,\epsilon}})$ is the class of the admissible functions, i.e., the set of the functions in $H^{1,2}(B,\R^3)\cap C^0(\partial B, \R^3)$ which map $\partial B$ weakly monotonic onto $\Gamma$, satisfy a three point condition and have an image almost everywhere in $\overline{\mathfrak{S}_{\beta+\delta,\epsilon}}$ (see also \cite{MinSurf}).

Since $Q$ verifies \eqref{eq1normlinftyQ} and $ \overline{\mathfrak{S}_{\beta+\delta, \e}} \subset \overline{\mathfrak{C}_{\beta+\delta}}$ we get that $\mathcal{F}$ is coercive. In fact, considering the associated Lagrangian $$e(p,q)=\frac{1}{2}(q_1^2+q_2^2) +2Q(p) \cdot q_1\wedge q_2,$$
where $p=(x,y,z) \in \overline{\mathfrak{S}_{\beta+\delta, \e}}$, $q=(q_1, q_2) \in \R^3\times\R^3$, by elementary computations and using \eqref{eq1normlinftyQ}, we get that, for any $p \in\overline{\mathfrak{S}_{\beta + \delta,\epsilon}}$, $q \in  \R^3\times\R^3$
$$\left(\frac{1}{2} - \frac{c_{\beta - \delta}}{2}\right) (|p_1|^2+|p_2|^2) \leq e(p,q) \leq \left(\frac{1}{2} + \frac{c_{\beta-\delta}}{2}\right)  (|p_1|^2+|p_2|^2).$$

In order to minimize the energy functional we have to prove that the class of admissible functions is not empty, i.e., $\mathcal{C}(\Gamma,\overline{\mathfrak{S}_{\beta+\delta,\epsilon}}) \neq \varnothing$.
To this end we recall that since $\Gamma$ is rectifiable it is well known that the set  $\mathcal{C}(\Gamma,\R^3)$ is not empty (see  \cite{MinSurf} pag. 255) and there exists a minimal surface $Y \in \mathcal{C}(\Gamma,\R^3)$ spanning $\Gamma$. Since $Y \in C^0(\overline{B}, \R^3) \cap C^2({B}, \R^3)$ is harmonic, by the Convex hull theorem (see Theorem 1, Section 4.1 of \cite{RegMinSurf}) we have that $Y(\overline{B})$ is contained in the convex hull of $\Gamma$. In particular, being $\overline{\mathfrak{S}_{\beta+\delta,\epsilon}}$ convex we get that $Y(\overline{B}) \subset\overline{\mathfrak{S}_{\beta+\delta,\epsilon}}$. Hence $Y \in  \mathcal{C}(\Gamma,\overline{\mathfrak{S}_{\beta+\delta,\epsilon}})$.

 By Theorem 3 in Section 4.7 of \cite{RegMinSurf} we have that the variational problem $\mathcal{P}(\Gamma,\overline{\mathfrak{S}_{\beta + \delta,\epsilon}})$ has a weak solution $X \in \mathcal{C}(\Gamma, \overline{\mathfrak{S}_{\beta + \delta,\epsilon}})$ and it satisfies the conformality relations
$$|X_u|^2=|X_v|^2, \ \ \ X_u \cdot X_v=0 \ \ \hbox{a.e. in}\ B. $$

\noindent\textbf{Step 6:} The weak solution $X$ found at Step 5 is a classical solution of \eqref{eqprob} and maps homeomorphically $\partial B$ onto $\Gamma$.\\ 
Since $\overline{\mathfrak{S}_{\beta + \delta,\epsilon}}$ is a closed and convex set, such that $\displaystyle\overline{\mathring{{\mathfrak{S}_{\beta + \delta,\epsilon}}}}={\mathfrak{S}_{\beta + \delta,\epsilon}}$, we have that $\overline{\mathfrak{S}_{\beta + \delta,\epsilon}}$ is a quasi-regular set (see Remark (i), pag 381 in \cite{RegMinSurf}). Thanks to a well known regularity result (see Theorem 4, pag 381 in \cite{RegMinSurf}) since $\overline{\mathfrak{S}_{\beta + \delta,\epsilon}}$ is quasi-regular it follows that $X$ is continuous up to the boundary.
In order to get more regularity and prove that $X$ is a classical solution of \eqref{eqprob}, we show that $X$ does not touch the boundary $\partial \mathfrak{S}_{\beta + \delta,\epsilon}=S_{\beta+\delta, \e}$. To prove this, we will argue by contradiction and use an important result of Gulliver and Spruck, which is a sort of geometric maximum principle. 

Assume by contradiction that $X$ touches $S_{\beta + \delta,\epsilon}$. The idea is to show that in this case we can construct $Y \in \mathcal{C}(\Gamma, \overline{\mathfrak{S}_{\beta + \delta,\epsilon}})$ such that $\mathcal{F}(Y)<\mathcal{F}(X)$, and hence, being $X$ of least energy we get a contradiction. To this end we define a ``truncation" map $T\colon\overline{\mathfrak{S}_{\beta + \delta,\epsilon}} \to \R^3$. 

In order to define $T$ we need some preliminary definitions: for $p \in \overline{\mathfrak{S}_{\beta + \delta,\epsilon}} $ we define $r(p):=dist(p,{S}_{\beta + \delta,\epsilon})$, we observe that there exists a neighborhood $V$ of ${S}_{\beta + \delta,\epsilon}$ such that for $p \in V$ there is a unique point $\pi(p) \in {S}_{\beta + \delta,\epsilon}$ with $|p-\pi(p)|=r(p)$. We observe that, in the definition of $\pi$ it is fundamental that ${S}_{\beta + \delta,\epsilon}$ is smooth: in fact, in the case of a cone, for any neighborhood $V$ of the cone, we have that any point $p \in V$ lying on the axis of the cone we have that $p$ is equidistant from ${S}_{\beta + \delta,\epsilon}$, so $\pi(p)$ cannot be defined as in the previous way. 

We also observe that $\pi:V \to {S}_{\beta + \delta,\epsilon}$ is a $C^1$ map. Finally, for $R>0$ we define $T\colon\overline{\mathfrak{S}_{\beta + \delta,\epsilon}} \to \R^3$ by setting $$T(p):= \begin{cases} \pi(p) + R N(\pi(p)) & \hbox{if} \ p \in V \ \hbox{and} \ r(p)\leq R,\\ p & \hbox{otherwise},\end{cases}$$ where $N(q)$ is the inward normal at $q \in {S}_{\beta + \delta,\epsilon} $. In general $T$ may be not continuous, but thanks to Theorem 3.1 of \cite{GulliverSpruck2}, since \eqref{eqstimaH} holds,  there exists $R_0 >0$ such that if $0<R\leq R_0$ and $\inf_{z \in B} r(X(w))< R$ we have $T \circ X \in C^0(B) \cap H^{1,2}(B,\R^3)$ and $\mathcal{F}(T\circ X)<\mathcal{F}(X)$. Since we are assuming that $X$ touches $S_{\beta + \delta,\epsilon}$ we have that $\inf_{z \in B} r(X(w))< R$, for any $0<R<R_0$. Hence $\mathcal{F}(T\circ X)<\mathcal{F}(X)$. It remains to prove that $T\circ X \in \mathcal{C}(\Gamma, \overline{\mathfrak{S}_{\beta + \delta,\epsilon}})$. From the proof of Theorem 3.1 in \cite{GulliverSpruck2} we know that $T \circ X \in C^0(B) \cap H^{1,2}(B,\R^3)$, moreover since $\Gamma$ is strictly in the interior of $ \overline{\mathfrak{S}_{\beta + \delta,\epsilon}}$, for $R$ sufficiently small, by definition of $T$ we have that $T(p)=p$, for any $p \in \Gamma$. Hence, since $X$ is a weakly monotonic map of $\partial B$ onto $\Gamma$, and satisfies a three point condition, the same holds for $T\circ X$, and thus $T\circ X \in \mathcal{C}(\Gamma, \overline{\mathfrak{S}_{\beta + \delta,\epsilon}})$ and we get the contradiction.

Therefore we have that $X(\overline B) \cap S_{\beta + \delta,\epsilon} = \varnothing$, so from Theorem 7 in Section 4.7 of \cite{RegMinSurf} we have that $X$ is a classical solution of \eqref{eqprob} and $X(\overline B) \subset \mathfrak{S}_{\beta+\delta, \e}$. Moreover we observe that by construction, for all sufficiently small $\e>0$ we have that ${\mathfrak{S}_{\beta + \delta,\epsilon}} \subset {\mathfrak{C}_{\beta + \delta}}$. 

We observe that $X\colon\partial B \to \Gamma$ is a homeomorphism. This follows in a standard manner (see for instance the proof of Theorem 8 and the Remark at page 402 in \cite{RegMinSurf}).\\

\noindent\textbf{Step 7:} $X(B) \subset {\mathfrak{C}_\beta}$.\\

We begin with proving that $X(\overline B) \subset \overline{{\mathfrak{C}_{\beta}}}\setminus \{0\}$. Let us set 
$\phi(u,v):=X(u,v) \cdot e_3 - |X(u,v)|\cos \beta$, where $e_3=(0, 0, 1)$. We want to prove that $\phi \geq 0$ in $\overline B$. To this end we first show that $-\Delta \phi \geq 0$ in $B$,  i.e., $\phi$ is super-harmonic in $B$.

By elementary computations we have that $\phi_u=X_u \cdot e_3 - \frac{X\cdot X_u}{|X|}\cos \beta$, and $$\phi_{uu}=X_{uu}\cdot e_3 - \frac{X_u^2+ X\cdot X_{uu}}{|X|} \cos \beta+\frac{( X\cdot X_{u})^2}{|X|^3} \cos \beta.$$
Hence we get that $$- \Delta \phi = -\Delta X\cdot e_3 + \frac{2E+ X\cdot \Delta X}{|X|} \cos \beta-\frac{( X\cdot X_{u})^2 + ( X\cdot X_{u})^2}{|X|^3} \cos \beta, $$
where $E=|X_u|^2=|X_v|^2$. Since $0 \not\in X(\overline B)$ we have that $\phi \in C^2(B)\cap C^0(\overline B)$. 
Now, recalling that $X$ is an H-surface, we deduce that in the subset $B^\prime \subset B$ of regular points it holds
\beq\label{eqphisuparm}
\begin{array}{lll}
- \Delta \phi &=& \displaystyle -\Delta X\cdot e_3 + \frac{2E+ X\cdot \Delta X}{|X|} \cos \beta-\frac{( X\cdot X_{u})^2 + ( X\cdot X_{v})^2}{|X|^3} \cos \beta\\[8pt]
&=& \displaystyle- 2H(X) (X_u \wedge X_v\cdot e_3) + \frac{2E}{|X|}\cos \beta+  2 H(X)(P(X)\cdot X_u \wedge X_v) \cos \beta\\[12pt]
&&   \displaystyle-\frac{( P(X)\cdot P(X_{u}))^2 + ( P(X)\cdot P(X_{v}))^2}{|X|}E \cos \beta\\[12pt]
&\geq& \displaystyle- 2|H(X)| E + \frac{2E}{|X|}\cos \beta-  2 |H(X)| E \cos \beta-\frac{E}{|X|} \cos \beta\\[12pt]
&=& \displaystyle  \frac{\cos \beta - 2(1+\cos \beta) |H(X)||X|}{|X|}E \geq 0~\!.
\end{array}
\eeq
We point out that the last inequality holds because $H$ satisfies assumption  \eqref{growthcondition}, and the previous one is a consequence of $(P(X)\cdot P (X_{u}))^2 + ( P(X)\cdot P(X_{u}))^2\leq 1$, which comes from the orthogonality of the versors $P(X_{u})$, $P(X_{v})$.

On the other hand, if $(u_0,v_0)\in B$ is a branch point then, from the first line of \eqref{eqphisuparm} we get that $-\Delta \phi (u_0,v_0)=0$. Hence we have proved that $-\Delta \phi \geq 0$ in $B$ and we are done.

Now, since $X$ maps $\partial B$ onto $\Gamma$ and $\Gamma \subset \overline{\mathfrak{C}_\beta} \setminus \{0\}$ we have that $\phi\geq 0$ on $\partial B$. Therefore, by the maximum principle we get that $\phi\geq 0$ in $\overline B$,  from which we get that $X(\overline B)\subset \overline {\mathfrak{S}_{\beta,\e}} \subset \overline{\mathfrak{C}_\beta}\setminus \{0\}$. 

Now, from Enclosure Theorem I (see Section 4.2 in \cite{RegMinSurf}), in view of \eqref{eqstimaH} (which holds for $\delta=0$), we get that $X(B) \subset \mathfrak{S}_{\beta, \e}$, from which we deduce that $X$ maps $B$ into $\mathfrak{C}_\beta$, and we are done.\\

 %
\noindent\textbf{Step 8:} Proof of \eqref{eqstimaH}.\\
Using the parametrization $\phi(t,\theta)=(\alpha_1(t)\cos \theta, \alpha_1(t) \sin \theta, \alpha_2(t))$ 
 we have that the mean curvature (with respect to the inward normal) of $S_{\beta + \delta,\epsilon}$ is given by 
$$ H_{S_{\beta + \delta,\epsilon}}=\frac{\alpha_1 (\alpha_1^{\prime}  \alpha_2^{\prime\prime}-\alpha_2^{\prime} \alpha_1^{\prime\prime}) + \alpha_2^{\prime} ((\alpha_1^{\prime})^2 + (\alpha_2^{\prime})^2)}{2\alpha_1 ((\alpha_1^{\prime})^2 + (\alpha_2^{\prime})^2)^{3/2}}$$
(see for instance \cite{AbateTovena}). For $t>t_\e$ we have that the rotation of $\sigma_\e(t)=(\alpha_1(t),0,\alpha_2(t))$ describes a portion of the cone $\partial\mathfrak{C}_{\beta+\delta}$ and it is elementary to see that
$$H_{S_{\beta + \delta,\epsilon}}(t)=\frac{\mathrm{cotan}(\beta+\delta)}{2t}.$$
On the other hand, if $p=\phi(t,\theta) \in S_{\beta + \delta,\epsilon}$, 
for $t>t_\e$, $\theta \in [0,2\pi[$ we have
$$|H(p)|\leq \frac{c_{\beta-\delta}}{|\phi(t,\theta)|}= \frac{c_{\beta-\delta}}{\sqrt{\alpha_1^2(t)+\alpha_2^2(t)}}= \frac{c_{\beta-\delta}}{t}.$$
Hence 
\beq\label{eq1prelstimaH}
|H(\phi(t,\theta))| < H_{S_{\beta + \delta,\epsilon}}(t),\ \hbox{for any} \ t>t_\e, \ \theta \in [0,2\pi[
\eeq
 if and only if
$$ c_{\beta-\delta} < \frac{\mathrm{cotan}{(\beta+\delta})}{2},$$
which holds true as displayed in Step 1.

For the remaining interval  $[0,t_\e]$ we have the following:
\beq\label{eqlimHetozero}
\lim_{\e \to 0} \min_{t \in [0,t_\e]} H_{S_{\beta + \delta,\epsilon}}(t) = + \infty.
\eeq

Before proving \eqref{eqlimHetozero} we observe that it implies that there exists a small $\bar\e>0$ such that $\|H\|_{\infty, \overline{\mathfrak{C}_{\beta + \delta}}\cap\{z\leq1\}} <  \min_{[0,t_\e]} H_{S_{\beta + \delta,\epsilon}}(t)$ for all $0<\e<\bar \e$. Hence for all sufficiently small $\e>0$ we have
\beq\label{eq2prelstimaH}
 H(\phi(t,\theta)) < H_{S_{\beta + \delta,\epsilon}}(t),\ \hbox{for any}\ t \in [0,t_\e],\ \theta \in [0,2\pi[.
\eeq
At the end, thanks to \eqref{eq1prelstimaH}, \eqref{eq2prelstimaH}, we get \eqref{eqstimaH}.

Now we prove \eqref{eqlimHetozero}. First, for $\e>0$ sufficiently small, for any $t \in [0,t_\e]$, we have that
\begin{equation}
\label{375}
H_{S_{\beta + \delta,\epsilon}}(t) \geq \frac{- 4\sqrt{3} \left(\frac{3}{8}\right)^4 \frac{ \cos^4(\beta+\delta)}{\epsilon^3} t ^2 + 4\sqrt{3} \left(\frac{3}{8}\right)^2  \frac{ \cos^2(\beta+\delta)}{\epsilon}  }{2\sin(\beta+\delta)  \left((\sin(\beta+\delta))^2 + (- 4\sqrt{3} \left(\frac{3}{8}\right)^4 \frac{ \cos^4(\beta+\delta)}{\epsilon^3} t ^3 + 4\sqrt{3} \left(\frac{3}{8}\right)^2  \frac{ \cos^2(\beta+\delta)}{\epsilon} t)^2\right)^{1/2}}.
\end{equation}
In fact observe that $\alpha_1^{\prime\prime}\equiv 0$, $\alpha_1\geq 0$, $\alpha_1^\prime >0$, $\alpha_2^{\prime\prime} > 0$ in $[0,t_\e]$, hence
\begin{eqnarray*}
H_{S_{\beta + \delta,\epsilon}}(t)&=&\frac{\alpha_1 (\alpha_1^{\prime}  \alpha_2^{\prime\prime}-\alpha_2^{\prime} \alpha_1^{\prime\prime}) + \alpha_2^{\prime} ((\alpha_1^{\prime})^2 + (\alpha_2^{\prime})^2)}{2\alpha_1 ((\alpha_1^{\prime})^2 + (\alpha_2^{\prime})^2)^{3/2}}\\[8pt]
 &\geq& \frac{\alpha_1 (\alpha_1^{\prime}  \alpha_2^{\prime\prime}) + \alpha_2^{\prime} ((\alpha_1^{\prime})^2 + (\alpha_2^{\prime})^2)}{2\alpha_1 ((\alpha_1^{\prime})^2 + (\alpha_2^{\prime})^2)^{3/2}}\geq \frac{ \alpha_2^{\prime} }{2\alpha_1 ((\alpha_1^{\prime})^2 + (\alpha_2^{\prime})^2)^{1/2}}
\end{eqnarray*}
that is \eqref{375}.
Now, setting $s:=\frac{t}{\e}$, $g_\e(s):= - 4\sqrt{3} \left(\frac{3}{8}\right)^4 \cos^4(\beta+\delta) \frac{s^2}{\e} + 4\sqrt{3} \left(\frac{3}{8}\right)^2  \cos^2(\beta+\delta) \frac{1}{\e} $, $h(s):=- 4\sqrt{3} \left(\frac{3}{8}\right)^4 \cos^4(\beta+\delta)  s^3 + 4\sqrt{3} \left(\frac{3}{8}\right)^2  \cos^2(\beta+\delta)s$ by the previous estimate we deduce that
$$ \min_{t \in [0,t_\e]}H_{S_{\beta + \delta,\epsilon}}(t) \geq \min_{s \in \left[0, \frac{8}{3\sqrt{3}}\frac{1}{\cos(\beta+\delta)}\right]}\frac{g_\e(s)}{2\sin(\beta+\delta)  \left((\sin(\beta+\delta))^2 + (h(s))^2\right)^{1/2}}.$$
Since $s \in \left[0, \frac{8}{3\sqrt{3}}\frac{1}{\cos(\beta+\delta)}\right]$, $g_\e(s)= 4\sqrt{3} \left(\frac{3}{8}\right)^2 \cos^2(\beta+\delta) \left(- \left(\frac{3}{8}\right)^2 \cos^2(\beta+\delta) \frac{s^2}{\e} + \frac{1}{\e}\right)$, it is elementary to see that $$g_\e(s) \geq  4\sqrt{3} \left(\frac{3}{8}\right)^2 \cos^2(\beta+\delta)\frac{2}{3\e}.$$
Moreover, since $s \mapsto\frac{1}{2\sin(\beta+\delta)  \left((\sin(\beta+\delta))^2 + (h(s))^2\right)^{1/2}}$  does not depend on $\e$, there exists a positive constant $C_1$ depending only on $\beta+\delta$ such that, for any $s \in \left[0, \frac{8}{3\sqrt{3}}\frac{1}{\cos(\beta+\delta)}\right]$, we have 
$$ \frac{1}{2\sin(\beta+\delta)  \left((\sin(\beta+\delta))^2 + (h(s))^2\right)^{1/2}}>C_1. $$
Finally, putting together these estimates, we have $$ \min_{t \in [0,t_\e]}H_{S_{\beta + \delta,\epsilon}}(t) \geq 4C_1\sqrt{3} \left(\frac{3}{8}\right)^2 \cos^2(\beta+\delta)\frac{2}{3\e} \to + \infty, \ \ \hbox{as} \ \e\to 0.$$
Hence \eqref{eqlimHetozero} is proved.\\
 
The proof is now complete.


\section{On $\beta$-convex domains and related results}
In this section we introduce the definition of $\beta$-convexity and prove some geometric results about $\beta$-convex subsets of $\mathbb{S}^2$ as well as geometric results about $H$-surfaces having support in a cone, whose boundary datum radially projects onto the boundary of a smooth $\beta$-convex subset.\\

Let $\Omega$ be a open subset of the unit sphere $\mathbb{S}^2$ such that $\partial \Omega$ is a Jordan curve.
We denote by $\mathfrak{C}_\Omega$ the conic region in $\R^3$ spanned by $\Omega$. 
\begin{definition}
We say that $\Omega$ is convex if $\mathfrak{C}_\Omega$ is a convex subset of $\R^3$.
\end{definition}

In order to get our results we need of a stronger convexity notion. For $\hat p_0 \in \mathbb{S}^2$ and $\beta \in (0,\frac{\pi}{2})$ we set 
$$\mathfrak{C}_{\hat p_0, \beta}:=\{x \in \R^3; \ x\cdot \hat p_0 - |x| \cos \beta >0\}.$$ We introduce the following definition:

\begin{definition}\label{defbetaconvx}
 Let $\beta \in (0,\frac{\pi}{2})$. We say that $\Omega$ verifies a $\beta$-cone condition at a given $p \in \partial \Omega$ if there exists $\hat p_0 \in \mathbb{S}^2$ such that $p \in \partial {\mathfrak{C}_{\hat p_0, \beta}}$ and $\overline\Omega \subset \overline{\mathfrak{C}_{\hat p_0, \beta}}$. We say that $\Omega$ is $\beta$-convex if, for any $p \in \partial \Omega$, $\Omega$ verifies a $\beta$-cone condition at $p$.
 \end{definition}

We observe that, by definition, if $\Omega$ is $\beta$-convex, then, it is strictly contained in a hemisphere.\\

At first sight, one could think that for any $p \in \partial \Omega$ there could be many $\hat p_0 \in \mathbb{S}^2$ satisfying the $\beta$-cone condition at $p$, but this is not the case:

\begin{proposition} \label{propclassadmvers}
Assume that $\Omega$ verifies a $\beta$-cone condition at $p \in \partial \Omega$ and that $\partial \Omega$ is a regular Jordan curve of class $C^1$. Then there exists only one $\hat p_0 \in \mathbb{S}^2$ such that $p \in \partial {\mathfrak{C}_{\hat p_0, \beta}}$ and $\overline\Omega \subset \overline{\mathfrak{C}_{\hat p_0, \beta}}$. Moreover the mapping $p \mapsto \hat p_0$ is continuous from $\partial \Omega$ into $\mathbb{S}^2$.
\end{proposition}

\begin{proof}
Let $\sigma:(-\delta,\delta) \to \partial \Omega$ be a $C^1$-parametrization of a portion of $\partial \Omega$, centered at $p$. Since $\partial \Omega$ is a regular curve we can assume that $\sigma^\prime(t) \neq 0$ in $(-\delta,\delta)$. Since $\Omega$ verifies a $\beta$-cone condition at $p$, then, all possible $\hat p_0=\hat p_0(p,\beta)$ lie in $\partial\mathfrak{C}_{p, \beta}\cap \mathbb{S}^2$. Now observe that for any admissible $\hat p_0$, since $\sigma(0) \cdot \hat p_0= \cos \beta$, $\sigma(t) \cdot \hat p_0 \geq \cos \beta$ in $(-\delta,\delta)$, then, the function $h(t):=\sigma(t) \cdot \hat p_0$ must have null derivative at $0$. Hence $\sigma^\prime(0) \cdot \hat p_0=0$, which means that all possible $\hat p_0(p,\beta)$ must lie in the plane $\{\sigma^\prime(0)\}^\perp$. We also observe that since $|\sigma|\equiv 1$, then, by deriving this relation, we get that $p \in \{\sigma^\prime(0)\}^\perp$.

Thus all possible $\hat p_0$ are given by the intersection $\partial\mathfrak{C}_{p, \beta} \cap \{\sigma^\prime(0)\}^\perp \cap \mathbb{S}^2$ which consists of two vectors $\hat p_{0,1}$, $\hat p_{0,2}$. By construction we observe that they generate two cones $\partial\mathfrak{C}_{\hat p_{0,1}, \beta}$, $\partial\mathfrak{C}_{\hat p_{0,2}, \beta}$ such that $\partial\mathfrak{C}_{\hat p_{0,1}, \beta} \cap \partial\mathfrak{C}_{\hat p_{0,2}, \beta}=\{\lambda p,\ \lambda \in \R^+ \}$. Hence, since $\overline\Omega$ must be entirely contained in one of the regions $\overline{\mathfrak{C}_{p_{0,1}, \beta}}$, $\overline{\mathfrak{C}_{p_{0,2}, \beta}}$, we have that only one of the two vectors $\hat p_{0,1}$, $\hat p_{0,2}$ is admissible. The first part of proof is then complete.

We prove now the continuity of  the map $p\mapsto \hat p_0$, from $\partial \Omega$ into $\mathbb{S}^2$.  If $\sigma:(-\delta,\delta)\to\partial \Omega$ is a local parametrization centered at $p \in \partial \Omega$, then, as seen in the first part of the proof we have $\hat p_0(\sigma(t))= \partial\mathfrak{C}_{\sigma(t), \beta} \cap \{\sigma^\prime(t)\}^\perp \cap \mathbb{S}^2\cap \overline \Omega$. Hence, it is clear that $\hat p_0(\sigma(t))$ depends continuously on $t$ and we are done.

The proof is then complete.
\end{proof}

Next proposition states that $\beta$-convexity is actually a convexity property. 
\begin{proposition}
If $\Omega$ is $\beta$-convex then $\Omega$ is convex.
\end{proposition}

\begin{proof}
Assume by contradiction that $\Omega$ is not convex. Then, there exist two distinct points $p_1, p_2 \in \mathfrak{C}_\Omega$ such that the segment $\sigma(t)$ joining $p_1$ and $p_2$ is not entirely contained in $\mathfrak{C}_\Omega$. Let us set $\hat p_1:=P(p_1)$, $\hat p_2:=P(p_2)$, $\hat \sigma:=P\circ\sigma$. Then, there exists $t_0 \in (0,1)$ such that $\hat \sigma (t_0) \in \partial \Omega$. Since $\Omega$ is $\beta$-convex, choosing $p:=\hat \sigma (t_0)$ in the definition, we get that there exists $\hat p_0$ such that $\overline \Omega$ is contained in the region $\overline{\mathfrak{C}_{\hat p_0, \beta}}$, and $p \in \partial \Omega \cap \partial\mathfrak{C}_{\hat p_0, \beta}$. We observe that since $\hat p_1, \hat p_2 \in \Omega$ then $\hat p_1, \hat p_2 \in \mathfrak{C}_{\hat p_0, \beta}$ (they cannot lie on its boundary  $\partial\mathfrak{C}_{\hat p_0, \beta}$, otherwise they would belong to $\partial \Omega$).

Hence we have that $p_1, p_2 \in \mathfrak{C}_{\hat p_0, \beta}$ but $\sigma (t_0) \not \in \mathfrak{C}_{\hat p_0, \beta}$ which contradicts the convexity of $\mathfrak{C}_{\hat p_0, \beta}$.
\end{proof}

Now let us examine the relationship between the notion of $\beta$-convexity and some geometrical properties of $H$-surfaces. We begin with the following preliminary result:

\begin{proposition}\label{propphipos}
Let $\beta \in (0,\frac{\pi}{2})$, let $\Omega \subset \mathbb{S}^2$ be a $\beta$-convex domain and let $\Gamma$ be a  smooth regular Jordan curve such that $P(\Gamma)\subset \partial\Omega$. Assume that $H$ satisfies  \eqref{growthcondition}, and let $X \in C^2(B, \R^3) \cap C^0(\overline{B},\R^3)$ be an H-surface, with $X(\overline{B})\subset \overline{\mathfrak{C}_{\beta}} \setminus\{0\}$. Then, for any $p \in \partial \Omega$, the associated function $\phi_{p}(u,v):=X(u,v) \cdot \hat p_0 - |X(u,v)|\cos \beta$ is strictly positive in $B$, where $\hat p_0=\hat p_0(p,\beta) \in \mathbb{S}^2$  is given by the definition of $\beta$-convexity. 
\end{proposition}

\begin{proof}
Let us fix $p \in \partial \Omega$. Since $\Omega$ is $\beta$-convex there exists $\hat p_0\in \mathbb{S}^2$ such that   $p \in \partial\mathfrak{C}_{\hat p_0, \beta}$ and $\overline\Omega$ is contained in $\overline{\mathfrak{C}_{\hat p_0, \beta}}$. Hence, setting $\phi_{p}(u,v):=X(u,v) \cdot \hat p_0 - |X(u,v)|\cos \beta$, we have that $\phi_{p}\geq 0$ in $\partial B$. By replacing $e_3$ with $\hat p_0$ in the proof of Step 7 we have that $\phi_p$ is super-harmonic in $B$, and by the maximum principle we get that $\phi_{p}\geq 0$ in $B$. From the strong maximum principle it follows that $\phi_{p}>0$ in $B$ or $\phi_{p}\equiv 0$ in $B$. To complete the proof we have to show that the latter possibility cannot occur.

Assume by contradiction that $\phi_{p}\equiv 0$ in $B$, then, by definition and since $X$ is smooth we have that $X(\overline{B}) \subset \partial\mathfrak{C}_{\hat p_0, \beta}\setminus\{0\}$. Without loss of generality we can assume that $\hat p_0=e_3$ so that $X(\overline{B})$ is entirely contained in the surface $\partial\mathfrak{C}_{\beta}\setminus\{0\}$ which is the surface of revolution generated by the rotation, with respect of the $z$-axis, of the curve $\sigma$, lying in the $xz$-plane, given by $\sigma(t):=(\alpha_1(t),0,\alpha_2(t))$, where $\alpha_1(t)=\sin (\beta) t$, $\alpha_2(t)=\cos (\beta) t$, $t>0$. 
As seen in the proof of Theorem \ref{thm:existence}, using the parametrization $\phi(t,\theta)=(\alpha_1(t)\cos \theta, \alpha_1(t) \sin \theta, \alpha_2(t))$, we have that the mean curvature of $\partial\mathfrak{C}_{\beta}\setminus\{0\}$ (with respect to the inward normal) is given by $H_{\partial\mathfrak{C}_{\beta}\setminus\{0\}}(t)=\frac{1}{t} \frac{\ \mathrm{cotan}(\beta)}{2}$, $t>0$, moreover $|H(\phi(t,\theta))| < H_{\partial\mathfrak{C}_{ \beta}\setminus\{0\}}(t)$ for all $t>0$, $\theta \in [0,2\pi]$. In fact, since $H$ satisfies \eqref{growthcondition}, then for all $p=\phi(t,\theta) \in  \partial\mathfrak{C}_{\beta}\setminus\{0\}$ we have 
\beq\label{eqconfrcurv}
|H(\phi(t,\theta))|\leq \frac{c_\beta}{|\phi(t,\theta)|} = \frac{c_\beta}{t}< \frac{1}{2t} \ \mathrm{cotan}(\beta)=H_{\partial\mathfrak{C}_{\beta}\setminus\{0\}}(t),
\eeq
because $c_\beta=\frac{\cos \beta}{2(1+\cos \beta)} < \frac{\mathrm{cotan} \beta}{2}$. Thanks to \eqref{eqconfrcurv}, Theorem 2 and Corollary 3 in Section 4.4 of \cite{RegMinSurf} (or by Enclosure Theorem I in Section 4.2 of \cite{RegMinSurf}) it follows that $X(B)\cap (\partial\mathfrak{C}_{\beta}\setminus\{0\}) = \varnothing$ which gives a contradiction. 
 The proof is then concluded.
\end{proof}

\begin{corollary}\label{cor1hopf}
Under the same assumptions of the previous proposition we have that, for any $(u,v) \in \partial B$, the normal derivative, with respect to the exterior normal $\nu$ of the function $\phi_p$, corresponding to $p=PX(u,v) \in \partial \Omega$, is strictly negative at $(u,v)$, i.e. $$\frac{\partial}{\partial \nu}\phi_p(u,v)<0.$$ In particular, if $X \in C^1(\overline B, \R^3)$ then $X$ has no boundary branch points.
\end{corollary}
\begin{proof}
Let us fix $(u,v) \in \partial B$ and let $p=PX(u,v) \in\partial\Omega$. Consider the associated function $\phi_p$. As seen in the proof of Proposition \ref{propphipos} we have that $\phi_p$ is super-harmonic in $B$ and  $\phi_p>0$ in B. Hence, since $\phi_p(u,v)=0$, by Hopf's Lemma, we get that $\frac{\partial}{\partial \nu}\phi_p(u,v)<0$, where $\nu=(\nu^1, \nu^2)$ denotes the exterior normal at $(u,v) \in \partial B$. The first part is then proved.

For the second part we observe that since $\frac{\partial\phi_p}{\partial u} = X_u \cdot \hat p_0 - \frac{X \cdot X_u}{|X|} \cos \beta$, we have
$$\frac{\partial}{\partial \nu}\phi_p(u,v)= (X_u \cdot \hat p_0) \nu^1 + (X_v \cdot \hat p_0) \nu^2 -  \frac{(X \cdot X_u) \nu^1 + (X \cdot X_v) \nu^2 }{|X|} \cos \beta < 0.$$
Since $(u,v)$ is arbitrary we get that $X$ cannot have branch points on $\partial B$.
\end{proof}
Another important and immediate consequence of Proposition \ref{propphipos} is the following:

\begin{proposition}\label{propprojcontomega}
Under the same assumptions of Proposition \ref{propphipos} we have that $$PX(B) \subset \Omega.$$
In particular $X$ has support in the cone spanned by $\Omega$, i.e., $X(\overline B) \subset \overline{\mathfrak{C}_\Omega} \setminus \{0\}$.
\end{proposition}

\begin{proof}
Assume by contradiction that there exists some $(u_0,v_0) \in B$ such that $PX(u_0,v_0) \in \mathbb{S}^2\setminus \Omega$, then, necessarily, there exists $(u_1,v_1) \in B$ such that $X(u_1,v_1) \in \partial\mathfrak{C}_\Omega\setminus\{0\}$. 

In fact, on the contrary, we would have that $X(B) \cap ( \partial\mathfrak{C}_\Omega\setminus\{0\}) = \varnothing$, and hence we would have $$X(B)=(X(B)\cap \mathfrak{C}_\Omega) \cup [X(B)\cap (\R^3\setminus \overline {\mathfrak{C}_\Omega})].$$
Since we are assuming that $PX(u_0,v_0) \in \mathbb{S}^2\setminus \Omega$, we have that both the open sets in the right-hand side are nonempty, hence, since they are disjoint and $X(B)$ is connected we get a contradiction.

Hence there exists $(u_1,v_1) \in B$ such that $X(u_1,v_1) \in \partial\mathfrak{C}_\Omega\setminus\{0\}$, and taking $p_1=PX(u_1,v_1) \in \partial \Omega$, by the definition of $\beta$-convexity and applying Proposition \ref{propphipos} to the function $\phi_{p_1}$, we get a contradiction since $\phi_{p_1}(u_1,v_1)=0$.
\end{proof}

\section{Stable H-surfaces with one-one radial projection onto a $\beta$-convex subset}

In this section we analyze the geometrical properties of stable $H$-surfaces whose boundary is a Jordan curve $\Gamma$ that projects bijectively onto the boundary of a smooth $\beta$-convex domain $\Omega$ of the unit sphere $\mathbb{S}^2$. It will be understood, if not specified, that $\Gamma$ is a Jordan curve of class $C^{3,\alpha}$ and $H \in C^{1,\alpha}$, for some $\alpha \in (0,1)$, so that the solution found in Theorem \ref{thm:existence} is of class $C^{3,\alpha}(\overline{B}, \R^3)$ (see also Proposition \ref{regpropbhdata}).

We begin with a preliminary proposition:

\begin{proposition}\label{propifsegnposboundaryth}
Let $X$ be a stable H-surface of class $C^{3,\alpha}(\overline{B}, \R^3)$, with $H$ satisfying 
\eqref{monassumption}. 
Assume that $N\cdot X>0$ on $\partial B$, then $N\cdot X>0$ in $\overline{B}$.
\end{proposition}
\begin{proof}
Let us set $f:=N\cdot X$. By elementary computations we have $f_u=N_u \cdot X + N\cdot X_u=N_u \cdot X$, and thus $f_{uu}=N_{uu}\cdot X + N_u \cdot X_u$. Deriving the relation $N \cdot X_u \equiv 0$ we also get that $N_u \cdot X_u=-N\cdot X_{uu}$. Hence  $f_{uu}= N_{uu}\cdot X - N\cdot X_{uu}$ and thus $\Delta f= \Delta N \cdot X - N\cdot \Delta X$. Now, thanks to Theorem \ref{teoeqnorm}, in the subset $B^\prime \subset B$ of regular points, we get that
\begin{equation*}
\begin{array}{lll}
\Delta f + 2p f&=&\displaystyle  - 2E\nabla H(X) \cdot X - 2H(X)[N \cdot (X_u \wedge X_v)] \\[6pt] 
&=&\displaystyle  - 2E(\nabla H(X) \cdot X  + H(X)).
\end{array}
\end{equation*}
Since we are assuming \eqref{monassumption} we have $- 2E(\nabla H(X) \cdot X  + H(X))\leq 0$ in $B^\prime$.\\
On the other hand in the subset of branch points of $X$ we have $\Delta f + 2p f=\displaystyle  - 2E\nabla H(X) \cdot X - 2H(X)[N \cdot (X_u \wedge X_v)] =0$. Now applying Proposition \ref{propstabsegnpos} (we recall that $p \in C^{0,\alpha}(B)$) 
we get that $f>0$ in $\overline{B}$ and we are done.
\end{proof}

It remains to study the sign of $N \cdot X$ on the boundary $\partial B$. The next proposition ensures that $N\cdot X$ never vanishes on $\partial B$. 

\begin{proposition}\label{propcostsign}
Let $\Omega$ be a $\beta$-convex domain of class $C^{3,\alpha}$ and let $\Gamma$ be a Jordan curve of class $C^{3,\alpha}$ which radially projects onto $\partial \Omega$. Assume that $X$ is an H-surface of class $C^{3,\alpha}(\overline B, \R^3)$ with $H$ satisfying \eqref{growthcondition}. Then the function $N\cdot X$ never vanishes on $\partial B$, hence  $N\cdot X$ has a constant sign on $\partial B$.
\end{proposition}

\begin{proof}
Let us set $f:=N\cdot X$. Assume by contradiction that there exists $z_0 \in \partial B$ such that $f(z_0)=0$. In particular, since $X$ has no boundary branch points (see Corollary \ref{cor1hopf}) then we have $X(z_0) \cdot X_u(z_0)\wedge X_v(z_0)=0$. This means that $X(z_0) \in Span \{ X_u(z_0),\ X_v(z_0)\}:=\Pi$. Hence it follows that
$$(PX)_u(z_0) = \frac{X_u(z_0)}{|X(z_0)|} - \frac{X(z_0) \cdot X_u(z_0)}{|X(z_0)|^3} X(z_0)  \in \Pi,$$
and the same happens for $(PX)_v(z_0)$. Moreover, by deriving $|PX|\equiv 1$, we get that $PX \cdot (PX)_u \equiv 0$, $PX \cdot (PX)_v \equiv 0$. 

Let us set $v_1:=PX(z_0)$, $v_2:=(PX)_u(z_0)$, $v_3:=(PX)_v(z_0)$ and let us observe that $v_1,v_2, v_3 \in \Pi$ and $v_1 \cdot v_2=v_1\cdot v_3=0$. In particular, since $v_1\neq 0$ we deduce that 
\beq\label{eqwedprodv}
v_2 \wedge v_3 =0.
\eeq

Let $\phi=X \cdot \hat p_0 - |X| \cos \beta$ be the associated function to $v_1 \in \partial \Omega$. In particular we have $\phi(z_0)=0$, $\phi \geq 0$ in $\overline{B}$ and $\frac{\partial \phi}{\partial \nu}(z_0) <0$ (see Corollary \ref{cor1hopf}). Let us also introduce the function $\psi\colon\R \to \R$, defined by $\psi(\theta):=PX(\cos \theta, \sin \theta) \cdot \hat p_0 - \cos \beta$ and let $\theta_0 \in [0,2 \pi[$ such that $z_0=(\cos \theta_0, \sin \theta_0)$. Since $\psi\geq 0$, $\psi(\theta_0)=0$ and $\psi \in C^1(\R)$ we have that $\psi^\prime (\theta_0)=0$. This means that 
\beq\label{eq1propsc}
v_2\cdot \hat p_0 (-\nu^2) + v_3 \cdot \hat p_0 (\nu^1) =0, 
\eeq
where $\nu_1=\cos \theta_0$, $\nu_2=\sin \theta_0$. Moreover, since $\phi(z_0)=0$ we observe that $\frac{\partial \phi}{\partial \nu}(z_0) <0$ can be rewrited as
\beq\label{eq2propsc}
v_2\cdot \hat p_0 (\nu^1) + v_3\cdot \hat p_0 (\nu^2) < 0.
\eeq
We show that \eqref{eqwedprodv}, \eqref{eq1propsc} and \eqref{eq2propsc} lead to a contradiction.

If $v_2 \neq 0$ and $v_3 \neq 0$ then, setting $a:=v_2 \cdot \hat p_0$, $b:=v_3 \cdot \hat p_0$ we rewrite \eqref{eq1propsc}, \eqref{eq2propsc} as
\begin{equation*}
\begin{cases}
\begin{array}{lll}
\displaystyle  -\nu^2 a +  \nu^1 b &=& \displaystyle0\\
\displaystyle \ \ \nu^1 a +  \nu^2 b &=& \displaystyle -k,
\end{array}
\end{cases}
\end{equation*}
for some $k>0$. Then, by elementary computations it follows that $(a,b)= - k (\nu^1,\nu^2).$ Hence we have that
\begin{equation}\label{sistprojp0}
\begin{cases}
\begin{array}{lll}
\displaystyle  v_2 \cdot \hat p_0 &=& - k \nu^1,\\
\displaystyle  v_3 \cdot \hat p_0 &=& - k \nu^2.
\end{array}
\end{cases}
\end{equation}
On the other hand $v_2 \wedge v_3 =0$ implies that $v_2=\lambda v_3$, for some $\lambda \neq 0$, and hence from \eqref{sistprojp0} we have
\beq\label{eq3propcs}
 -k \nu^1 = \lambda v_3 \cdot \hat p_0 = \lambda (- k \nu^2).
\eeq
Remembering that \eqref{eqprob} is invariant under conformal transformations of the unit disk into itself, up to a rotation of angle $2\pi - \theta_0$, we can assume that $z_0=(1,0)$, in particular $\nu^1=1$, $\nu^2=0$ and so, since $k\neq 0$, we contradicts \eqref{eq3propcs}.

It remains to examine the case in which at least one between $v_2$, $v_3$ is zero. Assume by contradiction that $v_2=0$, then, thanks to \eqref{eq1propsc}, \eqref{eq2propsc} we get that
\begin{equation}\label{eq4sistcontr}
\begin{cases}
\begin{array}{lll}
\displaystyle   v_3\cdot \hat p_0 (\nu^1) &=& 0\\
\displaystyle  v_3\cdot \hat p_0 (\nu^2) &<& 0,
\end{array}
\end{cases}
\end{equation}
Up to a rotation we can assume that $\nu^1\neq 0$, $\nu^2\neq 0$, and hence \eqref{eq4sistcontr} gives a contradiction. The same argument shows that $v_3=0$ cannot happen. The proof is complete.
\end{proof}

It remains to prove that $N \cdot X>0$ on $\partial B$. To this end we we introduce the following definition:

\begin{definition}\label{definitposorjc}
Let $\Omega \subset \mathbb{S}^2$ be a $\beta$-convex domain, such that $\partial \Omega$ is a regular Jordan curve of class $C^{1}$, i.e., there exists a parametrization $\gamma\colon\partial B \to \partial \Omega$ of class $C^1$ which is a homeomorphism and satisfies $\gamma^\prime(z)\neq 0$ for all $z=e^{i \theta} \in \partial B$, where $\gamma^\prime(z)=\frac{d}{d\theta} \gamma(e^{i\theta})$. We say that $\partial \Omega$ is positively oriented by $\gamma$ if we have $(\gamma^\prime(z) \wedge \gamma(z)) \cdot \hat p_0(z) <0$, for all $z=e^{i\theta}$, $\theta \in [0,2\pi[$, where $\hat p_0(z)$ is the versor associated to $\gamma(z)$, given by the definition of $\beta$-convexity.
\end{definition}

\begin{remark}\label{remarkposorient}
We point out that the sign of $(\gamma^\prime(z) \wedge \gamma(z)) \cdot \hat p_0(z)$ is well defined since, as proved in Proposition \ref{propclassadmvers}, there is only one  $\hat p_0(z) \in \mathbb{S}^2$ satisfying the $\beta$-convexity condition at $\gamma(z) \in \partial \Omega$. Moreover, for any $z\in \partial B$, it cannot happen that $(\gamma^\prime(z) \wedge \gamma(z)) \cdot \hat p_0(z) =0$.  In fact, if we consider the scalar function $\theta \mapsto  h(\theta):=\gamma(e^{i \theta}) \cdot \hat p_0$, since $h$ has a minimum at $\theta_0$ corresponding to $z$, we get that $\gamma^\prime(z) \cdot \hat p_0=0$ and hence if by contradiction $\hat p_0(z) \in Span \{\gamma^\prime(z), \gamma(z)\}$, then $\hat p_0$ would be proportional to $\gamma(z)$ which is not possible. Hence we must have $Det [\gamma^\prime(z), \gamma(z), \hat p_0(z)] \neq 0$ on $\partial B$. Furthermore, thanks to the second part of Proposition \ref{propclassadmvers}, we deduce that $Det [\gamma^\prime(z), \gamma(z), \hat p_0(z)]$ is continuous on $\partial B$. Hence Definition \ref{definitposorjc} well defines an orientation on $\partial \Omega$. 
\end{remark}

Now we have all the instruments to state our assumption, which will be crucial for getting our next results. 

Given an $H$-surface $X$ of class $C^{3,\alpha}(\overline{B}, \R^3)$ spanning a regular Jordan curve $\Gamma$ of class $C^{3,\alpha}$ we introduce the following:\\[4pt]
\textbf{Assumption (I)}:
\begin{itemize}
\item[(i)] $\Gamma$ is a radial graph, i.e. there exists a domain $\Omega\subset\mathbb{S}^{2}$ and a map $g\colon\partial\Omega\to\R^{+}$ (with the same regularity of $\Gamma$) such that $\Gamma=\{g(p)p~|~p\in\partial\Omega\}$;
\item[(ii)] the domain $\Omega$ is $\beta$-convex;
\item[(iii)] the radial projection of $X|_{\partial B}$ induces a positive orientation on $\partial\Omega$.
\end{itemize} 
 


\begin{remark} \label{remarkposorientms}
We observe that, in our context, assumption (iii) makes sense. In fact, by definition of $H$-surface we have that $X|_{\partial B}\colon\partial B \to \Gamma$ is an homeomorphism and by Corollary \ref{cor1hopf} we know that $X$ has no boundary branch points.
\end{remark}

\begin{proposition}\label{propsegnposboundary}
Let $\Gamma$ be a regular Jordan curve of class $C^{3,\alpha}$ contained in $\overline{\mathfrak{C}_\beta} \setminus \{0\}$ and let $X \in C^{3,\alpha}(\overline{B}, \R^3)$ be an $H$-surface spanning $\Gamma$. Suppose that Assumption (I) is satisfied. Then $N \cdot X>0$ on $\partial B$.
\end{proposition}

\begin{proof} 
Let $z_0=(u_0,z_0) \in \partial B$ the point in which $|X|^2$ achieves its maximum and set $M_0:=\sup_{p \in \Gamma} |p|^2$. Let $\hat p_0 \in \mathbb{S}^2$ be the versor associated to $PX(z_0)$ by the definition of $\beta$-convexity. Up to a rotation of angle $\theta_0 \in [0,2\pi[$ we can assume that $z_0=(1,0)$. We point out that this does not change the induced orientation on $\partial \Omega$. Thanks to Corollary \ref{cor1hopf}, since $\nu=(1,0)$, we have that
\beq\label{eq1Hoplem}
X_u(z_0) \cdot \hat p_0 < \frac{X(z_0)\cdot X_u(z_0)}{|X(z_0)|}\cos \beta.
\eeq

On the other hand if we consider the map $\eta\colon\R\to\R$ given by $\eta(\theta):=|X(\cos \theta, \sin \theta)|^2$, since $\theta_0=0$ is a maximum point and $X$ is smooth up to the boundary, then $\psi^\prime(0)=0$ and hence we get that
\beq\label{eq1perp1}
 X(z_0) \cdot X_v(z_0) =0.
\eeq
Now consider the function $\psi\colon\R\to\R$, given by $\psi(\theta):=X(\cos \theta, \sin \theta)\cdot \hat p_0 - |X(\cos \theta, \sin \theta)| \cos \beta$. Since $\theta=0$ is a minimum point for $\psi$, and $X$ is smooth up to the boundary, we get that $\psi^\prime(0)=0$, and taking into account of \eqref{eq1perp1}, we deduce that
\beq\label{eq2perp2}
 X_v(z_0) \cdot \hat p_0=0.
\eeq

Equations \eqref{eq1perp1}, \eqref{eq2perp2} mean that $X_v(z_0)$ is orthogonal to both $X(z_0)$ and $\hat p_0$. Thus, for some $\lambda \in \R\setminus\{0\}$, it holds $X_v(z_0)=\lambda \hat p_0 \wedge X(z_0)$. Thanks to assumption (I), being $\frac{X_v(z_0)}{|X_v(z_0)|}$ the tangent versor to $\Gamma$ at $X(z_0)$ (we recall that, by Corollary \ref{cor1hopf}, $X$ has no boundary branch points) we have that $\lambda>0$, in particular $X_v(z_0)$  has the same direction and verse of $\hat p_0 \wedge X(z_0)$. To prove this, we first observe that thanks to the definition of $\beta$-convexity $\hat p_0$ and $X(z_0)$ must be linearly independent, so setting $\Pi:=Span \{\hat p_0, X(z_0)\}$ we have that $\Pi$ is a plane. Moreover, taking into account of Assumption (I) and Remark \ref{remarkposorient}, we have that $PX$ induces a positive orientation on $\partial \Omega$. Hence, by Definition \ref{definitposorjc}, since $(PX)^\prime(0)=\frac{X_v(z_0)}{|X(z_0)|}-\frac{X(z_0)\cdot X_v(z_0)}{|X(z_0)|^3}X(z_0)$, we must have 
$$\left(\frac{X_v(z_0)}{|X(z_0)|} \wedge X(z_0)\right) \cdot \hat p_0 =Det\left[\frac{X_v(z_0)}{|X(z_0)|}, X(z_0), \hat p_0\right]<0.$$
Hence, being $X_v(z_0)=\lambda\  \hat p_0 \wedge X(z_0)$, by the elementary properties of the determinant we get that $\lambda>0$.


Now let us consider the map $|X|^2\colon\overline{B} \to \R$. Since $X$ is an H-surface, and $H$ satisfies \eqref{growthcondition} we have that $|X|^2$ is subharmonic. In fact, by elementary computations, we have $(|X|^2)_u=2X\cdot X_u$, $(|X|^2)_{uu}=2 X_u\cdot X_u + 2 X \cdot X_{uu}$ and hence 

\begin{equation*}
\begin{array}{lll}
\displaystyle  -\Delta |X|^2 &=& \displaystyle - 4 E - 4 X \cdot H(X) (X_u \wedge X_v)\\[6pt]
&\leq& \displaystyle -4 E + 4 |X||H(X)| |X_u \wedge X_v|\\[6pt]
&\leq& \displaystyle -(4 - 4 c_\beta) E \leq 0.
\end{array}
\end{equation*}

In particular $|X|^2 - M_0$ is subharmonic and $|X|^2 - M_0 \leq 0$. Hence, by Hopf's lemma, since $|X|^2-M_0\not\equiv 0$ (otherwise $X$ would be a portion of a sphere, and hence $|H(X)|\equiv \frac{1}{\sqrt{M_0}}$, which contradicts  \eqref{growthcondition}), we get that
\beq\label{eq2Hoplem}
 X(z_0) \cdot X_u(z_0)>0~\!.
\eeq

Now, let us observe that by construction and since $X_u \cdot X_v\equiv 0$ we have $X_u(z_0) \in \Pi$. We want to understand where $X_u$ is located with respect to $\hat p_0$ and $X(z_0)$. By construction the two vectors $\hat p_0$ and $X(z_0)$ determine an angle of amplitude $\beta$. Let us denote by $R_1$ the angular region in $\Pi$ generated by $\hat p_0$ and $X(z_0)$, and by $R_2$ its complementary in $\Pi$.

We show that $X_u(z_0) \not\in R_1$. In fact if $X_u(z_0) \in R_1$, then, denoting by $\alpha \in ]0,\beta]$ the angle between $X_u(z_0)$ and $X(z_0)$ ($\alpha\neq 0$ in view of Proposition \ref{propcostsign}, or by \eqref{eq1Hoplem}) we have that $\beta-\alpha$ is the angle between $\hat p_0$ and $X_u(z_0)$. Then, by dividing by $|X_u(z_0)|$ each side of \eqref{eq1Hoplem}, we get that
\beq\label{eq1riletta}
\cos (\beta - \alpha) < \cos (\alpha) \cos (\beta),
\eeq
 and, by elementary trigonometry, we see that this last inequality is contradictory since both $\alpha$ and $\beta$ are in $]0,\pi/2[$.
 
 Hence, we have that $X_u(z_0) \in R_2$ and thanks to \eqref{eq2Hoplem} $X_u(z_0)$ must also lie in the half-plane $T:=\{p \in \Pi; \ p\cdot X(z_0) > 0\}$. Thus, $X_u(z_0) \in R_2 \cap T$, and let us consider the two subregions in which $R_2 \cap T$ splits: $R_{2,1}$, $R_{2,2}$. $R_{2,1}$ is defined as the subset of $R_2 \cap T$ such that $\hat p_0 \in \partial R_{2,1}$. Arguing as in the previous case we see that $X_u(z_0) \not\in R_{2,1}$. In fact, if $X_u(z_0) \in R_{2,1}$, denoting by $\alpha \in ]\beta,\pi/2[$ the angle between $X_u(z_0)$ and $X(z_0)$ we have that $\alpha-\beta$ is the angle between $\hat p_0$ and $X_u(z_0)$ and as before we have
 $$\cos (\alpha-\beta) < \cos (\alpha) \cos (\beta),$$ 
 which is contradictory.
 
 At the end the only possibility is $X_u(z_0) \in R_{2,2}$. 
Now, since $X_v(z_0)=\lambda\  \hat p_0 \wedge X(z_0)$, by the elementary properties of the determinant we get that 
 $$X(z_0) \cdot (X_u(z_0) \wedge X_v(z_0))=\lambda^3 Det\left[X_u,X,X\wedge \hat p_0\right]$$
and since $X_u \in R_{2,2}$ we have that $\{X_u,X,X\wedge \hat p_0\}$ is a positively oriented base of $\R^3$. Hence we get that $X(z_0) \cdot X_u(z_0) \wedge X_v(z_0)>0$. Now, thanks to Proposition \ref{propcostsign}, we have that the function $X\cdot X_u \wedge X_v$ has a constant sign on $\partial B$. Hence $N \cdot X>0$ on $\partial B$ and the proof is complete.
\end{proof}
 
From Proposition \ref{propifsegnposboundaryth} and Proposition \ref{propsegnposboundary} we finally get the following:

\begin{proposition} \label{propnormscalxpos}
Let $\Gamma$ be a regular Jordan curve of class $C^{3,\alpha}$ contained in $\overline{\mathfrak{C}_\beta}\setminus \{0\}$ and let $X \in C^{3,\alpha}(\overline{B}, \R^3)$ be an $H$-surface spanning $\Gamma$. Suppose that Assumption (I) is satisfied. Then $N \cdot X>0$ in $\overline B$.
\end{proposition}

\section{Global invertibility of the radial projection}
In this section we prove that under our assumptions the radial projection of an $H$-surface is a homeomorphism, in particular it can be represented as a radial graph. At the end of this section we prove Theorem \ref{thm:radial-graph}.

\begin{theorem}\label{mainteoglobinv}
Let $\Gamma$ be a regular Jordan curve of class $C^{3,\alpha}$ contained in $\overline{\mathfrak{C}_\beta}\setminus \{0\}$  and let $X$ be a stable H-surface of class $C^{3,\alpha}(\overline{B}, \R^3)$ spanning $\Gamma$, with $H$ satisfying \eqref{growthcondition}, \eqref{monassumption}. Suppose that Assumption (I) is satisfied. 
Then $PX\colon\overline{B} \to \overline{\Omega}$ is a homeomorphism. 
\end{theorem}

\begin{proof}
The idea is to apply a classical result of global invertibility (see Theorem \ref{globalinvtheo}). We divide the proof in four steps.\\

\textbf{Step 1:}  $PX$ is a surjective map from $\overline B$ to $\overline \Omega$.

Since $X$ maps homeomorphically $\partial B$ onto $\Gamma$, and $\Gamma$ satisfies assumption (I) then $PX$ maps homeomorphically $\partial B$ onto $\partial \Omega$ (it is a composition of a homeomorphism and a continuous  bijective map from a compact space into a Hausdorff space which is a homeomorphism too). 

 Without loss of generality, assume that $PX(\overline B)$ is contained in the upper hemisphere $\mathbb{S}^+:=\mathbb{S} \cap \{z>0\}$ and let us denote by $\pi\colon\mathbb{S}^2\setminus P_S \to \R^2$ the stereographic projection from the south pole $P_S=(0,0,-1)$. Since $\pi(PX)$ maps homeomorphically $\partial B$ onto $\pi(\partial \Omega)$ it follows that for $deg(\pi(PX),q)\equiv 1$ or $deg(\pi(PX),q)\equiv -1$, where $q \in \pi(PX(B))$. 
 
 In fact $q \not\in \pi( PX(\partial B))$ and by the basic properties of the degree (see for instance \cite{FonsecaGangbo}) we know that $q \mapsto deg(\pi(PX),q)$ is constant in each connected component of $\R^2\setminus \pi(PX(\partial B))$ (we recall that since $ \pi(PX(\partial B))$ is a Jordan curve then $\R^2\setminus \pi(PX(\partial B))$ has only two connected components), in particular it is constant for $q \in \pi(\Omega)$. Hence, being $\pi(PX)\colon\partial B \to \partial \Omega$ a homeomorphism there are only two possibilities: $deg(\pi(PX),q)\equiv 1$ or $deg(\pi(PX),q)\equiv -1$.
 
 Now we know that $PX \in C^0(\overline{B},\R^3)$, $P(B) \subset \Omega$ and being  $deg(\pi(PX),q)\neq 0$ for any  $q \in \pi(\Omega)$, it follows that $\pi(PX)(B)=\pi(\Omega)$ (see  \cite{FonsecaGangbo}). Being $\pi$ a diffeomorphism it follows that $PX(B)=\Omega$. At the end, since $PX$ maps $\partial B$ onto $\partial \Omega$, we have $PX(\overline B)=\overline \Omega$. Hence $PX$ is a surjective map from $\overline B$ to $\overline \Omega$.\\

\textbf{Step 2:}  $PX\colon\overline B \to \overline \Omega$ is locally invertible.

We begin with the local invertibility of $PX:B \to \Omega$. Being $(PX)_u= \frac{X_u}{|X|} - \frac{X\cdot X_u}{|X|^3}X$, by elementary computations we have 
\beq\label{eq1plinind}
(PX)_u\wedge (PX)_v \cdot PX = \frac{X_u \wedge X_v \cdot X}{|X|^3}.
\eeq
Hence, thanks to Proposition \ref{propnormscalxpos}, since $N \cdot X >0$ it follows that in the set $B^\prime$ of regular points of $X$ it holds
$$X_u \wedge X_v \cdot X>0, $$
which, in view of \eqref{eq1plinind}, implies that $(PX)_u(z)$ and $(PX)_v(z)$ are linearly independent in $B^\prime$. Thanks to a standard argument based on the the inverse function theorem it follows that $PX$ is a local diffeomorphism except on a discrete set of critical points (given by the branch points of $X$). Hence, from the standard properties of the degree (see for instance Theorem 2.9 in \cite{FonsecaGangbo}), from Proposition \ref{propnormscalxpos} and since $deg(\pi(PX),q)\equiv \pm 1$ for $q \in \pi(\Omega)$, it follows that each regular value has exactly one pre-image. In fact let $q \in \pi(\Omega)$ be a regular value, then the set of pre-images of $q$ is discrete and hence, being $\overline B$ compact, it is finite, and assuming for instance that $deg(\pi(PX),q)=1$ (see the proof of Step 1), by the index formula (see Theorem 2.9-(1) in  \cite{FonsecaGangbo})  we get that 
\beq\label{formulafonseca}
1=deg(\pi(PX),q)=\sum_{p \in [\pi(PX)]^{-1}(q)} i(\pi(PX),p).
\eeq
Now, being $q$ a regular value we have that $PX$ is local diffeomorphism at any $p \in [\pi(PX)]^{-1}(q)$, and $i(\pi(PX),p)=\pm 1$. Thanks to Proposition \ref{propnormscalxpos} and \eqref{eq1plinind} it follows that $(PX)_u\wedge (PX)_v \cdot PX>0$ in the set of regular points, in particular this holds near each $p \in [\pi(PX)]^{-1}(q)$. Hence, near each pre-image of $q$, $PX$ has the same orientation, so it follows that $i(\pi(PX),p)\equiv 1$ or  $i(\pi(PX),p)\equiv -1$, for  $p \in [\pi(PX)]^{-1}(q)$. Thus, from \eqref{formulafonseca}, we deduce that $i(\pi(PX),p)\equiv 1$ and 
$$1=deg(\pi(PX),q)=\sum_{p \in [\pi(PX)]^{-1}(q)} i(\pi(PX),p)=k,$$
where $k \in \mathbb{N}^+$ is the cardinality of the set $[\pi(PX)]^{-1}(q)$. Hence the only possibility is $k=1$, i.e.  $q$ has only one pre-image.
 
  It remains to prove the local invertibility in the finite set of branch points. Let $z_0$ be a branch point and assume by contradiction that $PX$ is not invertible at $z_0$ and set $p_0:=PX(z_0)$, then, for any neighborhood $V$ of $z_0$ we have that $PX$ is not injective in $V$. Since branch points are isolated we can assume without loss of generality that $V$ contains only $z_0$ as branch point. Then, there exist $z_1,z_2 \in B$, $z_1\neq z_2$ such that $PX(z_1)=PX(z_2)$ and necessarily one of them (for instance $z_1$) is not a branch point. It cannot happen that $PX(z_1)$ is a regular value, since we have proved that each regular value has exactly one pre-image. Hence $PX(z_1)=PX(z_0)=p_0$. By induction, repeating this argument we can construct a sequence of regular points $(z_{n}) \subset V$, $z_n \to z_0$ and such that $z_i\neq z_j$ for any $i\neq j$. In particular $S:=PX^{-1}(p_0)$ is not finite. Now, up to an isometry we can assume that $PX(z_0)=e_3$ and ${N}(z_0) \cdot e_3>0$. Let us denote by $\Pi_{e_3}\colon\R^3\to \R^2$ the projection of the first two coordinates. We observe that for any $z \in S$ we have $\Pi_{e_3}(X(z))=0$. On the other hand, arguing as in the proof of Theorem 1, Section 7.1 in \cite{MinSurf} (in particular, see (25)), since ${N}(z_0) \cdot e_3>0$, using the complex notation we can expand $\Pi_{e_3}(X(z))$ near $z_0$ as $\Pi_{e_3}(X(z))= l((z-z_0)^{n+1})+o(|z-z_0|^{n+1})$, where $n=n(z_0) \in \mathbb{N}$ is given by Theorem \ref{HWexpansion} and $l\colon\C \to \C$ is the map associated to a nonsingular real matrix (see (24), Sect. 7.1, \cite{MinSurf}). From this expansion we deduce that $0$ must have a finite set of pre-images near $z_0$, and hence we get a contradiction. Hence $PX:B \to \Omega$ is locally invertible. 
  
  Now we show that even $PX\colon\overline B \to \overline \Omega$ is locally invertible. In fact, as proved in Corollary \ref{cor1hopf}, $PX$ has no boundary branch points, so, considering a suitable $C^1$-extension of $PX$, to some open neighborhood $V$ of $\partial B$ we can assume that $X$ has no branch points in $V$. Now, from \eqref{eq1plinind}, Proposition \ref{propnormscalxpos}, we have
$$X_u \wedge X_v \cdot X>0 \ \ \hbox{in} \ V $$
which, in view of \eqref{eq1plinind}, implies that $(PX)_u(z)$ and $(PX)_v(z)$ are linearly independent in $V$. Hence, as before by an application of the inverse function theorem it follows that $PX$ is a locally invertible for any $z \in V$ and we are done.\\

\textbf{Step 3:} $PX\colon\overline B \to \overline \Omega$ is proper.

For any compact subset  $K \subset \overline \Omega$ we have that $K$ is closed and being $PX$ continuous we have $(PX)^{-1} (K)$ is a closed subset of $\overline{B}$. Being $\overline B$ compact it follows that $(PX)^{-1} (K)$ is compact.\\

\textbf{Step 4:}  $\overline\Omega$ is simply connected.

Thanks to an important result of differential geometry, known as Sch\"onflies's Theorem or also Jordan-Sch\"onflies's Theorem (for the proof see for instance \cite{Thomassen}) we know that the closure of the complement of the bounded region determined by a planar Jordan curve is homeomorphic to a closed ball, in particular it is simply connected. Hence, taking the stereographic projection of $\overline\Omega$, since $\partial \Omega$ is mapped onto a plane Jordan curve, we get that $\overline\Omega$ is simply connected.\\

From Step 1-Step 4, and being $\overline B$ arcwise connected we have that the hypotheses of Theorem \ref{globalinvtheo} are satisfied, so we get that $PX$ is a homeomorphism. The proof is complete.
\end{proof}

\begin{remark}
An immediate consequence of the previous theorem is that the H-surface $X$ can be expressed as a radial graph. In fact, let $(PX)^{-1}: \overline \Omega \to \overline B$ the inverse function of $PX$, and set $F(p):=(PX)^{-1}(p)$. Then, being $X(F(p))=PX(F(p)) |X(F(p))|= p |X(F(p))|$, we get that $$X(\overline B)=\{q \in \R^3;\ \ q= \lambda(p) p,\ p\ \in \overline\Omega\},$$ where $\lambda\colon\overline\Omega \to \R^+$  is the function defined by $\lambda(p):=|X(F(p))|$.\\
\end{remark}

\begin{proof}[Proof of Theorem \ref{thm:radial-graph}]
Let $X$ be the $H$-surface given by Theorem \ref{thm:existence}.  Under our assumptions we have that $X \in C^{3,\alpha}(\overline B, \R^3)$ (see Proposition \ref{regpropbhdata}) and $X$ is stable (see Remark \ref{stabilitymin}).
From a remarkable result of Gulliver (see Theorem 8.1 and Theorem 8.2 in \cite{Gulliver}) we have that $X$  is free of interior branch points and Corollary \ref{cor1hopf} excludes also boundary branch points. Hence, from the proof of Theorem \ref{mainteoglobinv} it follows that the radial projection of $X$ is actually a global diffeomorphism.
\end{proof}


\end{document}